\documentclass[12pt]{article}

\usepackage{geometry}
\usepackage[cp1251]{inputenc}
\usepackage[T2A]{fontenc}
\usepackage[english]{babel}

\usepackage{amsmath}
\usepackage{amssymb}
\usepackage{amsthm}
\usepackage{array}
\usepackage[auth-sc]{authblk}
\usepackage{bm}
\usepackage{cite}
\usepackage{color}
\usepackage{dsfont}
\usepackage{enumerate}
\usepackage{enumitem}
\usepackage{epsf}
\usepackage{euscript}
\usepackage{graphicx} % [dvips] "ломает" цвет текста в hyperref
\usepackage{indentfirst}
\usepackage{latexsym}
\usepackage{mathrsfs}
\usepackage{mathtools}
\usepackage{tikz}
\usepackage{pgf}

\newtheorem{theorem}             {Theorem}
\newtheorem{proposition}[theorem]{Proposition}

\newtheorem{corollary}  [theorem]{Corollary}

\newcommand{\logic}[1]{\mathbf{#1}}
\newcommand{\otuple}[1]{\langle{#1}\rangle}
\newcommand{\kframe}[1]{\mathfrak{#1}}
\newcommand{\kFrame}[1]{\boldsymbol{\kframe{#1}}}
\newcommand{\kmodel}[1]{\mathfrak{#1}}
\newcommand{\kModel}[1]{\boldsymbol{\kmodel{#1}}}
\newcommand{\nat}{\ensuremath{\mathds{N}}}

\begin{document}

%%%

\title{Algorithmic complexity of \protect\\ monadic multimodal predicate logics with equality \protect\\ over finite Kripke frames\thanks{The work on the paper was partially supported by the HSE Academic Fund Programme, Project~\mbox{23-00-022}.}}

\author[1]{I.~Agadzhanian}
\author[2]{M.~Rybakov}
\author[3]{D.~Shkatov}

\affil[1]{HSE University}
\affil[2]{IITP RAS, HSE University, Tver State University}
\affil[3]{University of the Witwatersrand, Johannesburg}

\date{}

%%%
\maketitle

\thispagestyle{empty}

\section{Introduction}

Monadic modal and superintuitionistic logics are, as a rule,
undecidable in very poor vocabularies---in most cases, to prove
undecidability, it suffices to use a single monadic predicate letter
and two or three individual
variables~\cite{RSh19SL,RSh20AiML,RShJLC20a,RShJLC21b,RShJLC21c,RybIGPL22}
(for undecidability of related fragments of the classical logics,
see~\cite{TG87,RSh19SAICSIT,MR:2022DM,MR:2023LI}). At the same time, the monadic
fragment with equality of the classical predicate logic $\logic{QCl}^=$ is
decidable~\cite{BGG97}.
%both in languages with at most two individual
%variables (and any set of predicate letters) and languages with
%monadic predicate letters and equality only (and any set of individual
%variables)
Hence, it is of interest to identify settings where decidability can
be obtained. 

Proofs of undecidability of monadic fragments usually rely on the
so-called the ``Kripke trick''~\cite{Kripke62}, a simulation of a
subformula $P(x,y)$ of a classical formula with a monadic modal
formula $\Diamond(Q_1(x)\wedge Q_2(y))$.  Hence, indentifying
decidable fragments involves discovering setting where the Kripke
trick is not applicable.  This has been done syntactically by Wolter
and Zakharyaschev~\cite{WZ01}, who discovered monodic fragments
(note that these differ from monadic fragments) disallowing
the application of modalities to formulas with more than one
parameter.  Here, we consider a simple semantical setting where the
Kripke trick does not work: the monadic predicate logic with equality
of a Kripke frame with finitely many possible worlds (but, possibly,
infinite domains).  We also obtain precise complexity bounds for
monadic logics of classes of Kripke frames with finitely many possible
worlds. This is of interest since precise bounds beyond $\Sigma^0_1$
hardness are scarce in the literature on predicate modal logic.  The
observations presented here are generalizations to the multimodal
settings of results from in~\cite{MR:2017LI,RSh:2023JLC}.

\section{Preliminaries}

We consider the $n$-modal, where $n\in\mathds{N}^+$, predicate
language $\mathcal{L}_n$ obtained by adding to the classical predicate
language $\mathcal{L}_0$ unary modal connectives
$\Box_1,\ldots,\Box_n$, as well as the language $\mathcal{L}^=_n$
obtained by adding to $\mathcal{L}_n$ a designated binary predicate
letter $=$. The definitions of formulas are standard. A
\textit{monadic $\mathcal{L}_n$-formula} contains only monadic
predicate letters.  A \textit{monadic $\mathcal{L}_n$-formula with
  equality} contains only monadic predicate letters and $=$.

By a normal $n$-modal predicate logic we mean a set of
$\mathcal{L}_n$-formulas including the classical predicate logic
$\mathbf{QCl}$ and the minimal normal $n$-modal propositional logic
$\mathbf{K}_n$ and closed under Modus Ponens, Substitution,
Necessitation, and Generalisation.  A normal $n$-modal predicate logic
with equality additionally contains the classical equality axioms.
The minimal logic containing $\mathbf{QCl}$ and the $n$-modal
propositional logic $L$ is denoted by $\mathbf{Q}L$; the minimal
extension of $\mathbf{Q}L$ containing the classical equality axioms is
denoted by $\mathbf{Q}^=L$.  The minimal extension of an $n$-modal
predicate logic $L$ containing, for each $k \in \{1, \ldots, n\}$, the
Barcan formula
$\bm{bf}_k = \forall x\, \Box_k P(x) \to \Box_k \forall x\, P(x)$, is
denoted by $L.\mathbf{bf}$.

%The modality of 1-modal logics, propositional or predicate, is denoted
%by $\Box$.
A \textit{fusion} of $1$-modal propositional logics
$L_1, \ldots, L_n$ is the logic
$L_1 \ast \ldots \ast L_n = \mathbf{K}_n \oplus (L_1\cup L'_2 \cup \ldots
\cup L'_n)$, where $L'_i$ is obtained from $L_i$ by replacing every
occurrence of $\Box_1$ with~$\Box_i$.

We use the framework of Kripke semantics for logics with and without
equality (for more details, see~\cite{GShS}; our terminology differs from
that adopted in~\cite{GShS}).  There are two natural way to extend the
well-known Kripke semantics for logics without equality to logics with equality;
to treat equality as identity or as hereditary congruence.  Unlike the
classical logic, these two treatments of equality are not equivalent:
the formula $x\ne y\to \Box_k(x\ne y)$ is valid if $=$ is interpreted
as identity, but not valid if $=$ is interpreted as hereditary
congruence.  Here, except in Section~\ref{sec:discussion}, we treat
equality as congruence.

A Kripke $n$-frame is a tuple $\kframe{F}=\otuple{W,R_1,\ldots,R_n}$,
where $W$ is a non-empty set of worlds and $R_1,\ldots,R_n$
are binary accessibility relations on~$W$. An augmented
$n$-frame is a tuple $\kFrame{F} = \otuple{\kframe{F}, D}$, where
$\kframe{F}$ is a Kripke $n$-frame and $D$ a family $(D_w)_{w \in W}$
of non-empty domains satisfying the expanding domains
  condition: for every $w, v \in W$,
\begin{itemize}
\item[]
$
\begin{array}{lrcl}
  (\mathit{E})~ & wR_k v & \Longrightarrow & D_w\subseteq D_v. \\
\end{array}
$
\end{itemize}
The condition ($\mathit{E}$) is required for soundness and
completeness of predicate modal logics whose
\mbox{$\mathcal{L}_0$-fragment} is $\mathbf{QCl}$.
If an augmented $n$-frame satisfies
\begin{itemize}
\item[]
$
\begin{array}{lrcl}
  (\mathit{C})~ & wR_k v & \Longrightarrow & D_w = D_v, \\
\end{array}
$
\end{itemize}
then it is called a locally constant augmented $n$-frame.  A model is a
tuple $\kModel{M}=\otuple{\kFrame{F},I}$, where $\kFrame{F}$ is an
augmented $n$-frame and $I$ is a family $(I_w)_{w\in W}$ of
interpretations of predicate letters: $I_w(P)\subseteq D_w^m$, for
every $m$-ary letter~$P$.

An augmented $n$-frame with equality is a tuple
$\kFrame{F} = \otuple{\kframe{F}, D, {\equiv}}$, where
$\langle \kframe{F}, D \rangle$ is an augmented $n$-frame and
${\equiv}$ is a family $(\equiv_w)_{w \in W}$ of equivalence
relations, with ${\equiv_w} \subseteq D_w^2$ whenever $w \in W$,
satisfying the heredity condition: for every $w, v \in W$,
\begin{itemize}
\item[]
$
\begin{array}{lrcl}
  (\mathit{H})~ & wR_k v & \Longrightarrow & {\equiv_w} \subseteq {\equiv_v}.
\end{array}
$
\end{itemize}
The condition $(\mathit{H})$ corresponds to the formula
$x=y\to \Box_k(x=y)$, which belongs to $\mathbf{Q}^=\logic{K}$, and hence to
every normal modal predicate logic with equality. A~model with
equality is a tuple $\kModel{M}=\otuple{\kFrame{F},I}$, where
$\kFrame{F}$ is an augmented $n$-frame with equality and $I$ is a
family $(I_w)_{w\in W}$ of interpretations of predicate letters such
that $\equiv_w$ is a congruence on the classical model
$M_w=\otuple{D_w,I_w}$.

The truth relation for $\mathcal{L}_n$ and $\mathcal{L}^=_n$ is defined by usual way; in
particular, if $a,b\in D_w$ and $\bar{c}$ is a list of elements of $D_w$ of a suitable length, then
$$
\begin{array}{lcl}
  \kModel{M},w\models a = b 
    & \leftrightharpoons 
    & a \equiv_w b; \\
  \kModel{M},w\models P(\bar{c}) 
    & \leftrightharpoons 
    & \bar{c}\in I_w(P); 
    \\
  \kModel{M},w\models \forall x\,\varphi(x,\bar{c}) 
    & \leftrightharpoons 
    & \mbox{$\kModel{M},w\models \varphi(d,\bar{c})$, for every $d\in D_w$;} \\
  \kModel{M},w\models \Box_k\varphi(\bar{c}) 
    & \leftrightharpoons 
    & \mbox{$\kModel{M},v\models \varphi(\bar{c})$, for every $v\in R_k(w)$.} \\
\end{array}
$$
The following definitions concern both $\mathcal{L}_n$ and
$\mathcal{L}^=_n$; for the latter, all the models and augmented frames
should be understood as those with equality.  A formula $\varphi$ is
true at a world $w$ if a universal closure of $\varphi$ is true
at~$w$. A formula~$\varphi$ is true in a model~$\kModel{M}$ if
$\varphi$ true at every world of~$\kModel{M}$; $\varphi$~is valid on
an augmented $n$-frame~$\kFrame{F}$ if it is true in every model
over~$\kFrame{F}$; $\varphi$~is valid on a Kripke
$n$-frame~$\kframe{F}$ if $\varphi$ is valid on every augmented
$n$-frame over~$\kframe{F}$; $\varphi$~is valid on a
class~$\mathscr{C}$ of augmented frames if it is valid on every
augmented frame from $\mathscr{C}$.

If~$\mathscr{C}$ is a class of Kripke $n$-frames and $\kframe{F}$ is a
Kripke frame, then
\begin{itemize}
\item $L(\mathscr{C})$ denotes the set of $\mathcal{L}_n$-formulas
  valid on~$\mathscr{C}$;
\item $L_c(\mathscr{C})$ denotes the set of $\mathcal{L}_n$-formulas
  valid on every locally constant augmented $n$-frame over a Kripke
  frame from~$\mathscr{C}$;
\item $L^=(\mathscr{C})$ denotes the set of $\mathcal{L}^=_n$-formulas
  valid on~$\mathscr{C}$;
\item $L^=_c(\mathscr{C})$ denotes the set of
  $\mathcal{L}^=_n$-formulas valid on every locally constant augmented
  $n$-frame with equality over a Kripke frame from~$\mathscr{C}$.
\end{itemize}
We write $L^=(\kframe{F})$ and $L_c^=(\kframe{F})$ rather than $L^=(\{\kframe{F}\})$ and $L_c^=(\{\kframe{F}\})$, respectively.

A Kripke $n$-frame $\langle W, R_1, \ldots, R_n \rangle$ is
finite if $W$ is a finite set.  If $L$ is an $n$-modal
predicate logic (with or without equality), then $L^{\mathit{wfin}}$ denotes the set of formulas
valid on every finite Kripke frame validating~$L$; this set is a
normal $n$-modal predicate logic.

\section{Main results}

The following is our main technical result:

%\begin{lemma}
%If a monadic modal formula with equality containing $k$ monadic
%predicate letters and $m$ individual variables is refuted on an
%augmented frame with equality
%$\boldsymbol{\mathfrak{F}} = \langle W, R_1, \ldots, R_n, D,
%\{\equiv_w\}_{w \in W} \rangle$ with finite $W$, then it is refuted on
%an augmented frame with equality
%${\boldsymbol{\mathfrak{F}}}' = \langle W, R_1, \ldots, R_n, D',
%\{\equiv'_{w}\}_{\mathstrut{w \in W}} \rangle$ with
%$|{D'}^+| \leqslant |W| \cdot m\cdot 2^{|W| (k+1)}$; moreover, if
%${\boldsymbol{\mathfrak{F}}}$ is locally constant, then so is
%${\boldsymbol{\mathfrak{F}}}'$.% and
%\end{lemma}

%Hence, we immediately obtain the following:

\begin{proposition}
  \label{thr:finite-class-eq}
  Let $\kframe{F}$ be a finite Kripke frame.  Then, the monadic fragments
  with equality of the logics $L^=(\kframe{F})$ and
  $L^=_c(\kframe{F})$ are both decidable.
\end{proposition}

From Proposition~\ref{thr:finite-class-eq} we obtain the following:

\begin{theorem}
  \label{thr:Pi^0_1-membership-eq}
  Let $\mathscr{C}$ be a recursively enumerable class of finite Kripke
  $n$-frames.  Then the monadic fragments with equality of the logics
  $L^=(\mathscr{C})$ and $L^=_c(\mathscr{C})$ are both in $\Pi^0_1$.
\end{theorem}

It is known~\cite[Theorem 3.9]{RShJLC20a} that, if $L$ is a logic from
one of the intervals
$[\mathbf{QK}^{\mathit{wfin}}, \mathbf{QGL.3.bf}^{\mathit{wfin}}]$,
$[\mathbf{QK}^{\mathit{wfin}}, \mathbf{QGrz.3.bf}^{\mathit{wfin}}]$ or
$[\mathbf{QK}^{\mathit{wfin}}, \mathbf{QS5}^{\mathit{wfin}}]$, then
the monadic fragment of $L$ is $\Pi^0_1$-hard.  Hence,
Theorem~\ref{thr:Pi^0_1-membership-eq}, together with the transfer of
completeness theorem for fusions~\cite[Theorem 4.1]{GKWZ}, give us the
following:

\begin{corollary}
  \label{thr:Pi^0_1-eq-fusions}
  Let $L = L_1 \ast \ldots \ast L_n$, where $L_1, \ldots, L_n$ are normal $1$-modal propositional logics,  such that
  \begin{itemize}
  \item $L_i \subseteq \mathbf{S5}$ or
    $L_i \subseteq \mathbf{GL.3}$ or $L_i \subseteq \mathbf{Grz.3}$, for some $i \in \{1, \ldots, n \}$;
  \item the class of finite Kripke frames validating $L$ is recursively
    enumerable.
  \end{itemize}
  Then, the monadic fragments
  of the logics\/ $\mathbf{Q}L^{\mathit{wfin}}$,
  $\mathbf{Q}L.\mathbf{bf}^{\mathit{wfin}}$ and the monadic fragments with equality of the logics
  $\mathbf{Q}^=L^{\mathit{wfin}}$,\/
  $\mathbf{Q}^=L.\mathbf{bf}^{\mathit{wfin}}$ are all\/
  \mbox{$\Pi^0_1$-complete}.
\end{corollary}

\begin{corollary}
  \label{thr:Pi^0_1-eq}
  Let $L$ be one of the logics\/ $\mathbf{K}$, $\mathbf{T}$,
  $\mathbf{D}$, $\mathbf{K4}$, $\mathbf{K4.3}$, $\mathbf{S4}$,
  $\mathbf{S4.3}$, $\mathbf{GL}$, $\mathbf{GL.3}$, $\mathbf{Grz}$,
  $\mathbf{Grz.3}$, $\mathbf{KB}$, $\mathbf{KTB}$, $\mathbf{K5}$,
  $\mathbf{K45}$, $\mathbf{S5}$.  Then, the monadic fragments of the logics\/ $\mathbf{Q}L_n^{\mathit{wfin}}$,
  $\mathbf{Q}L_n.\mathbf{bf}^{\mathit{wfin}}$ and the monadic fragments with equality of the logics
  $\mathbf{Q}^=L_n^{\mathit{wfin}}$,\/
  $\mathbf{Q}^=L_n.\mathbf{bf}^{\mathit{wfin}}$ are\/
  \mbox{$\Pi^0_1$-complete}.
\end{corollary}

%We notice that all logics with $\Pi^0_1$-complete monadic fragments
%mentioned in
%Corollaries~\ref{thr:Pi^0_1-eq-fusions}--\ref{thr:Pi^0_1-eq-products}
%are $\Sigma^0_1$-hard in languages with binary
%letters~\cite{RShJLC20a}.

From Proposition~\ref{thr:finite-class-eq} we also obtain the
following:

\begin{theorem}
  \label{thr:fixed-branching}
  Let $\mathscr{C}$ be a decidable class of Kripke $1$-frames closed
  under the operation of taking subframes and satisfying the condition
  that there exists $m \in \nat$ such that $|R(w)| \leqslant m$
  whenever $\langle W, R \rangle \in \mathscr{C}$ and $w \in W$.  Then
  the monadic fragments of the logics $L(\mathscr{C})$,
  $L_c(\mathscr{C})$ and the monadic fragments with equality of the logics $L^=(\mathscr{C})$, $L^=_c(\mathscr{C})$ are decidable.
\end{theorem}

Recall that $\logic{Alt}_n$ is a monomodal logic complete with respect
to the class of Kripke frames where every world sees at most
$n$~worlds. Using completeness of the predicate counterpart of
$\logic{Alt}_n$~\cite{ShehtmanShkkatov20}, which, using~\cite[Theorem
3.8.7]{GShS}, implies the completeness of $\logic{Q}^=\logic{Alt}_n$, we obtain
the following:
%statement
%for predicate counterparts of $\logic{Alt}_n$.

\begin{theorem}
  \label{thr:alt-eq}
  The monadic fragments of logics\/ $\mathbf{QAlt}_n$,
  $\mathbf{QAlt}_n\mathbf{.bf}$, $\mathbf{Q}^=\logic{Alt}_n$, and
  $\mathbf{Q}^=\logic{Alt}_n\mathbf{.bf}$ are all decidable.
\end{theorem}

%Notice that all the statements remain true if we remove the equality from the languages.

\section{Discussion}
\label{sec:discussion}

Proposition~\ref{thr:finite-class-eq} and
Theorem~\ref{thr:Pi^0_1-membership-eq} remain true in the Kripke
semantics with equality as identity.  Proposition~\ref{thr:finite-class-eq}
and Theorem~\ref{thr:Pi^0_1-membership-eq} can be extended to logics
of frames with distinguished worlds.  All the results remain true if
logics $\logic{Q}L.\mathbf{bf}$ and $\logic{Q}^=L.\mathbf{bf}$ are replaced with logics in which, for some $k$, the Barcan formula $\bm{bf}_k$ is replaced with $\Box\bm{bf}_k$, where $\Box$ is a finite sequence of $\Box_1, \ldots, \Box_n$.  Lastly, we
note that similar results can be obtained for superintuitionistic
monadic logics~\cite{RSh:2023JLC}.

%

%\medskip
%\paragraph{\itshape Acknowledgements.} %optional
%The work on the paper was partially supported by the HSE Academic Fund Programme (Project~\mbox{23-00-022}).

%%%

%\bibliographystyle{plain}
%\bibliography{sources} 

\begin{thebibliography}{10}

\bibitem{BGG97}
Egon B\"{o}rger, Erich Gr\"{a}del, and Yuri Gurevich.
\newblock {\em The Classical Decision Problem}.
\newblock Springer, 1997.

\bibitem{GKWZ}
Dov Gabbay, Agi Kurucz, Frank Wolter and Michael Zakharyaschev.
\newblock {\em Many-Dimensional Modal Logics}, volume 48 of
  {\em Studies in Logic and the Foundations of Mathematics}.
\newblock Elsevier, 2003.

%\bibitem{GSh93}
%Dov Gabbay and Valentin Shehtman.
%\newblock Undecidability of modal and intermediate first-order logics with two
%  individual variables.
%\newblock {\em The Journal of Symbolic Logic}, 58(3):800--823, 1993.

\bibitem{GShS}
Dov Gabbay, Valentin Shehtman, and Dmitrij Skvortsov.
\newblock {\em Quantification in Nonclassical Logic, Volume 1}, volume 153 of
  {\em Studies in Logic and the Foundations of Mathematics}.
\newblock Elsevier, 2009.

%\bibitem{KKZ05}
%Roman Kontchakov, Agi Kurucz, and Michael Zakharyaschev.
%\newblock Undecidability of first-order intuitionistic and modal logics with
%  two variables.
%\newblock {\em Bulletin of Symbolic Logic}, 11(3):428--438, 2005.

\bibitem{Kripke62}
Saul~A. Kripke.
\newblock The undecidability of monadic modal quantification theory.
\newblock {\em Zeitschrift f\"{u}r Matematische Logik und Grundlagen der
  Mathematik}, 8:113--116, 1962.

\bibitem{MR:2017LI}
Mikhail Rybakov.
\newblock Undecidability of modal logics of unary predicate.  
\newblock \emph{Logical Investigations},
23(2):60--75, 2017. (In Russian)

\bibitem{MR:2022DM}
Mikhail Rybakov.
\newblock Computational complexity of binary predicate theories with a small number of variables in the language.
\newblock \emph{Doklady Mathematics},
\newblock 507(6):61--65, 2022.

\bibitem{MR:2023LI}
Mikhail Rybakov.
\newblock Binary predicate, transitive closure, two-three variables: shall we play dominoes?
\newblock \emph{Logical Investigations}, 29(1):114--146, 2023. (In Russian)

\bibitem{RybIGPL22}
Mikhail Rybakov.
\newblock Predicate counterparts of modal logics of provability: High
  undecidability and {K}ripke incompleteness.
\newblock To appear in \textit{{L}ogic {J}ournal of the {I}{G}{P}{L}}.

\bibitem{RSh19SL}
Mikhail Rybakov and Dmitry Shkatov.
\newblock Undecidability of first-order modal and intuitionistic logics with
  two variables and one monadic predicate letter.
\newblock {\em Studia Logica}, 107(4):695--717, 2019.

\bibitem{RSh19SAICSIT}
Mikhail Rybakov and Dmitry Shkatov.
\newblock Trakhtenbrot theorem for classical languages with three individual variables. 
\newblock \emph{Proceedings of the South African Institute of Computer Scientists and Information Technologists 2019}, Skukuza, South Africa, September 17--18, 2019, ACM, NY, USA, Article No.\,19:1--7. 2019.

\bibitem{RSh20AiML}
Mikhail Rybakov and Dmitry Shkatov.
\newblock Algorithmic properties of first-order modal logics of the natural
  number line in restricted languages.
\newblock In Nicola Olivetti, Rineke Verbrugge, Sara Negri, and Gabriel Sandu,
  editors, {\em Advances in Modal Logic}, volume~13. College Publications,
  2020.

\bibitem{RShJLC20a}
Mikhail Rybakov and Dmitry Shkatov.
\newblock Algorithmic properties of first-order modal logics of finite {K}ripke
  frames in restricted languages.
\newblock {\em Journal of Logic and Computation}, 30(7):1305--1329, 2020.


\bibitem{RShJLC21b}
Mikhail Rybakov and Dmitry Shkatov.
\newblock Algorithmic properties of first-order superintuitionistic logics of
  finite {K}ripke frames in restricted languages.
\newblock {\em Journal of Logic and Computation}, 31(2):494--522, 2021.

\bibitem{RShJLC21c}
Mikhail Rybakov and Dmitry Shkatov.
\newblock Algorithmic properties of first-order modal logics of linear {K}ripke
  frames in restricted languages.
\newblock {\em Journal of Logic and Computation}, 31(5):1266--1288, 2021.


\bibitem{RSh:2023JLC}
Mikhail Rybakov and Dmitry Shkatov.
\newblock Algorithmic properties of modal and
superintuitionistic logics of monadic predicates
over finite frames. 
\newblock \emph{Submitted to the Journal of Logic and Computation.}


\bibitem{ShehtmanShkkatov20}
Valentin Shehtman and Dmitry Shkatov.
\newblock Some prospects for semiproducts and products of modal logics.
\newblock In N.~Olivetti and R.~Verbrugge, editors, {\em Short Papers Advances
  in Modal Logic {A}i{M}{L} 2020}, pages 107--111. University of Helsinki,
  2020.

\bibitem{TG87}
Alfred Tarski and Steven Givant.
\newblock {\em A Formalization of Set Theory without Variables}, volume~41 of
  {\em American Mathematical Society Colloquium Publications}.
\newblock American Mathematical Society, 1987.

\bibitem{WZ01}
Frank Wolter and Michael Zakharyaschev.
\newblock Decidable fragments of first-order modal logics.
\newblock {\em The Journal of Symbolic Logic}, 66:1415--1438, 2001.

\end{thebibliography}

%%%

\end{document}